\documentclass[aip,reprint,jcp]{revtex4-1}

\usepackage{amsmath}
\usepackage{amsfonts}
\usepackage{amssymb}
\usepackage{amsthm}
\usepackage{chemarrow}
\usepackage{graphicx}
\usepackage{epstopdf}
\usepackage{algorithm}
\usepackage{algpseudocode}
\usepackage{bm}
\usepackage{adjustbox}
\usepackage{multirow}

\newcommand{\R}{\mathbb R}
\newcommand{\mbf}[1]{\mathbf{#1}}
\newcommand{\trans}{^\text{T}}

\theoremstyle{definition}

\draft 

\begin{document}


\title{Dynamic Bounds on Stochastic Chemical Kinetic Systems Using Semidefinite Programming} 



\author{Garrett R. Dowdy}
\author{Paul I. Barton}
\email[]{pib@mit.edu}
\affiliation{Department of Chemical Engineering, Massachusetts Institute of Technology, Cambridge, MA 02139, USA}


\date{\today}

\begin{abstract}
Applying the method of moments to the chemical master equation (CME) appearing in stochastic chemical kinetics often leads to the so-called closure problem.
Recently, several authors showed that this problem can be partially overcome using moment-based semidefinite programs (SDPs).
In particular, they showed that moment-based SDPs can be used to calculate rigorous bounds on various descriptions of the stochastic chemical kinetic system\rq{}s stationary distribution(s) -- for example, mean molecular counts, variances in these counts, and so on.
In this paper, we show that these ideas can be extended to the corresponding dynamic problem, calculating time-varying bounds on the same descriptions.


\end{abstract}


\maketitle 



%
%

%

\section{Introduction}
A stochastic chemical kinetic system is inherently uncertain.
Thus, rather than talking about \emph{the} state of the system, it is more natural to talk about the probability of each reachable state.
Considering all of these probabilities collectively, we have a probability \emph{distribution} over the set of reachable states.
This probability distribution changes over time, and they way it changes is governed by the chemical master equation (CME).
Computing the solution to this equation would give us a complete dynamic description of the system, specifying the probability of each state throughout time.
However, for most systems of practical importance, direct numerical solution of the CME is difficult, because the number of equations and variables (i.e., states) is very large, even infinite. \cite{higham2008modeling}


The classical strategy for dealing with this problem of the large number of states is to sample the reaction system using Gillespie\rq{}s Stochastic Simulation Algorithm (SSA).
While this algorithm is intuitively appealing and very easy to implement, it is often too slow in practice  \cite{gillespie2007stochastic}.
Many variants of Gillespie\rq{}s algorithm have been developed with the aim of increasing its speed.
Most of these involve some approximation that renders their results inexact and potentially misleading.
Those that retain the exactness of Gillespie\rq{}s algorithm remain fundamentally limited in that they must simulate every reaction \cite{constantino2016modeling}.

Another strategy for dealing with the large number of states is to give up trying to calculate the time-varying probability associated with each state and focus instead on summary descriptions of the probability distribution -- for example, mean molecular counts and variances in these counts.
Conveniently, these quantities can be expressed in terms of the \emph{moments} of the distribution.
Furthermore, one can use the CME to derive an ordinary differential equation (ODE) describing how the moments of the system change over time.\cite{smadbeck2013closure, sotiropoulos2011analytical, gillespie2009moment}
Unfortunately, this ODE usually suffers from the so-called “closure problem”, in which the time evolution of the moments up to order $m$ depends on the values of moments up to order $m+1$.

To deal with the closure problem, various authors have proposed “closure scheme” approximations \cite{gillespie2009moment, smadbeck2013closure}.  
While these approximations have some intuitive appeal, they generally cannot provide bounds on the error they introduce.
One notable exception is the closure scheme described by Naghnaeian and Del Vecchio\cite{naghnaeian2017robust} which can provide error bounds under the condition that the molecular count of each species present in the system is bounded.
However, the scalability of this method is doubtful from a theoretical perspective, as it requires solving a linear program (LP) whose size is proportional to the number of reachable states.


Recently, several authors \cite{dowdy2017using, dowdy2018bounds, sakurai2017convex, kuntz2017rigorous, ghusinga2017exact} independently proposed an alternative to closure schemes, describing a method for calculating rigorous bounds on several quantities of interest for steady-state (i.e. stationary) stochastic chemical kinetic distributions.
The central idea of this method was to adapt Lasserre\rq{}s\cite{lasserre2010moments} moment-based semidefinite programs (SDPs) to the problem of stochastic chemical kinetics.
 
In the present paper, we will extend this idea to calculate  \emph{time-varying} bounds on dynamic stochastic chemical kinetic systems.
Once again, these bounds will be obtained by solving moment-based SDPs.

\section{Mathematical Background}


\subsection{Mathematical Notation}
Throughout this paper, the symbol $\mathbb N$ will be used to denote the set of natural numbers $\{0, 1, 2, \dots \}$, the symbol $\mathbb Z$ will be used to denote the integers $\{\dots, -2, -1, 0, 1, 2, \dots \}$, and  $\R$ will be used to denote the real numbers.
Bold symbols will be used to represent vectors and matrices.
The dimensions of these vectors and matrices will be specified as they are introduced.
The vector $\mbf e_i=(0, \dots ,1, \dots 0)$ is the $i$th coordinate vector, in which all components are zero, except the $i$th component, which is $1$.
Angular brackets \lq\lq{}$\langle \cdot \rangle$\rq\rq{} will be used to denote an \lq\lq{}expected value\rq\rq{} or mean of a random variable.
The meanings of all other symbols should be clear from the context.

\subsection{Stochastic Chemical Kinetics Notation}
Consider a stochastic chemical kinetic system with $N$ distinct chemical species and $R$ reactions.  
The state of the system at time $t$ is described by the random vector $\mbf X(t)=(X_1(t), \dots ,X_N(t) ) \in \mathbb N^N$, 
where $X_i(t) \geq 0$ is the count of molecules of species $i$ present.  

The state changes with the occurrence of each reaction.  
For example, if $\mbf s_r \in \mathbb Z^N$ is the vector of stoichiometric coefficients of reaction $r$, and the system is in state $\mbf x \in \mathbb N^N$, then an occurrence of reaction $r$ takes the system to state $\mbf x+ \mbf s_r$. 
By chaining together multiple reactions, a system initially in state $\mbf X(0) \equiv \mbf x_0 \in \mathbb N^N$ can reach many possible states -- sometimes infinitely many.
Let this set of \emph{reachable states} be denoted with the symbol $\mathcal X \subset \mathbb N^N$.
A generic element of this set will be denoted $\mbf x \in \mathcal X$.



\subsection{Invariants and Independent Species} \label{invariants}
The stoichiometry matrix for the system is constructed by bringing together the stoichiometry vectors: $\mbf S \equiv [ \mbf s_1 \dots \mbf s_R ] \in \mathbb Z^{N \times R}$.
Often, this matrix will have a nontrivial left null space.
Let $\{ \mbf b_1, \dots, \mbf b_L \} \subset \R^N$ be a basis for this left null space. 
It can be shown that each of these vectors corresponds to an \emph{invariant} of the reaction system\cite{fjeld1974reaction} -- i.e., some linear combination of molecular counts that is constant with time.
In particular,
\begin{equation} \label{definition of invariant values}
\mbf b_j\trans \mbf X(t) = f_j, \ \ \ \forall j \in \{1, \dots, L\}, \ \ \ \forall t \geq 0,
\end{equation}
where each $f_j \in \R$ is a constant which we will call the \emph{value} of the $j$th invariant.
In what follows, we will assume that these invariant values are known.
This is true, for example, if we know the initial state $\mbf X(0)$, because the invariant values can then be calculated via Equation \eqref{definition of invariant values}.
However, our method does not rely explicitly on knowledge of the initial state $\mbf X(0)$.
This has some interesting implications regarding uncertainty in the initial state, which will be explored further in Section \ref{uncertain initial state section}.

If we set $\mbf f \equiv (f_1, \dots, f_L)$, and 
\begin{equation}
\mbf B \equiv 
\left [
\begin{array}{c}
\mbf b_1\trans \\
\vdots \\
\mbf b_L\trans \\
\end{array}
\right ] \in \R^{L \times N},
\end{equation}
then Equation \eqref{definition of invariant values} can be expressed concisely as
\begin{equation} \label{concise expression of invariants}
\mbf B \mbf X(t) = \mbf f, \ \ \ \forall t \geq 0.
\end{equation}
These equations imply that the set of reachable states $\mathcal X$ is contained in an affine subspace, i.e., that  $\mathcal X \subset \{ \mbf x \in \R^n : \mbf B \mbf x = \mbf f \}$.
Furthermore, they imply that not all molecular counts $X_1, \dots, X_N$ can vary independently.
To see this, let $\tilde{\mbf B} \in \R^{L \times L}$ be a matrix obtained by concatenating $L$ linearly independent columns of $\mbf B$, and let $\tilde{\mbf X}(t) \in \mathbb N^L$ be the vector of the corresponding components of $\mbf X(t)$.
Similarly, let $\hat{\mbf B} \in \R^{L \times \hat N}$ be the matrix obtained by concatenating the remaining $N - L \equiv \hat N$ columns of $\mbf B$, and let $\hat{\mbf X}(t) \in \mathbb N^{\hat N}$ be the vector of the corresponding components of $\mbf X(t)$.
Then, Equation \eqref{concise expression of invariants} can be rewritten as 
\begin{equation} \label{expanded expression of invariants}
\tilde{\mbf B} \tilde{\mbf X}(t) + \hat{\mbf B} \hat{\mbf X}(t) = \mbf f, \ \ \ \forall t \geq 0.
\end{equation}
By construction, $\tilde{\mbf B}$ is invertible, so if $\hat{\mbf X}(t)$ is known, this equation can be solved for $\tilde{\mbf X}(t)$:
\begin{equation} \label{expanded expression of invariants}
\begin{aligned}
\tilde{\mbf X}(t) 
&= \tilde{\mbf B}^{-1} \left (\mbf f - \hat{\mbf B} \hat{\mbf X}(t) \right ), \ \ \ \forall t \geq 0, \\
&= \tilde{\mbf B}^{-1} \mbf f - \tilde{\mbf B}^{-1} \hat{\mbf B} \hat{\mbf X}(t), \ \ \ \forall t \geq 0. \\
\end{aligned}
\end{equation}
Thus, knowing $\hat{\mbf X}$ is enough to know the state of the system.
We can think of the chemical species whose molecular counts are specified in the vector $\hat{\mbf X}$ as being the \lq\lq{}independent species\rq\rq{}.
In general, there will be several possible ways to pick $L$ linearly independent columns of $\mbf B$. This means that we have some flexibility in choosing which species to treat as independent.

\subsection{A Reduced State Space}
Every full-dimensional reachable state $\mbf x \in \mathcal X \subset \mathbb N^N$ has a corresponding \emph{reduced reachable state}, $\hat{\mbf x} \in \mathbb N^{\hat N}$, obtained by selecting the counts of the independent species from $\mbf x$.
We will denote the set of all these reduced reachable states as $\hat{\mathcal X} \subset \mathbb N^{\hat N}$.
Similarly, for every stoichiometry vector $\mbf s_r \in \mathbb Z^N$, there is a corresponding reduced stoichiometry vector $\hat{\mbf s}_r \in \mathbb Z^{\hat N}$, obtained by selecting the components of $\hat{\mbf s}_r$ corresponding to the independent species.

Working in the reduced state space is computationally convenient because it focuses attention on the variables in the stochastic chemical kinetic system that are actually independent and can thus reduce the dimension of the problems we want to solve.
For the sake of brevity, in what follows, we will often loosely refer to the reduced state as simply the \lq\lq{}state\rq\rq{}.
That we are in fact referring to the reduced state should be clear from the context.

We know that the molecular counts of the independent species must be nonnegative.
So for any $\hat{\mbf x} \in \hat{\mathcal X}$, we must have $\hat{\mbf x} \geq \mbf 0$.
Furthermore, we know that the molecular counts of the dependent species must be nonnegative.
By Equation \eqref{expanded expression of invariants}, this implies $\tilde{\mbf B}^{-1} \mbf f - \tilde{\mbf B}^{-1} \hat{\mbf B} \hat{\mbf x} \geq \mbf 0$.
It follows that the set of reduced reachable states $\hat{\mathcal X}$ must be contained in the following polyhedral set:
\begin{equation} \label{polyhedral set}
\bar{\mathcal X} \equiv 
\left \{ \hat {\mbf x} \in \R^{\hat N} :
\begin{array}{c}
\hat{\mbf x} \geq \mbf 0, \\
\tilde{\mbf B}^{-1} \mbf f - \tilde{\mbf B}^{-1} \hat{\mbf B} \hat{\mbf x} \geq \mbf 0
\end{array} 
\right \}.
\end{equation}

\subsection{The Chemical Master Equation} \label{cme}
Because of the stochastic nature of the system, there is some uncertainty as to the (reduced) state at time $t$, and we express this uncertainty by assigning a probability $\text{Pr}(\hat{\mbf X}(t) = \hat{\mbf x},t) \equiv P(\hat{\mbf x},t)$ to each of the reachable states $\hat{\mbf x} \in \hat{\mathcal X}$.  
This probability distribution $P(\cdot,t)$ changes over time according to the chemical master equation (CME):
\begin{equation} \label{cme equation}
\begin{aligned}
\frac{dP}{dt} (\hat{\mbf x},t) =\sum_{r=1}^R [P(\hat{\mbf x} - \hat{\mbf s}_r,t)a_r (\hat{\mbf x} - &\hat{\mbf s}_r )-P( \hat{\mbf x},t) a_r (\hat{\mbf x})] , \\
&\ \ \ \forall \hat{\mbf x} \in \hat{\mathcal X},
\end{aligned}
\end{equation}
where $a_r$ is the \lq\lq{}propensity function\rq\rq{} of reaction $r$.
The details of this propensity function are described in Higham\cite{higham2008modeling}.
However, we want to point out two things:
first, $a_r(\cdot)$ is always a polynomial in $\hat{\mbf x}$;
second, $a_r$ is proportional to a rate constant $c_r$ for reaction $r$.
This $c_r$ is not necessarily the same as the macroscopic rate constant $k_r$ one would use in deterministic chemical kinetics, but there is a connection between the two constants. 
See Higham\cite{higham2008modeling} and Gillespie\cite{gillespie1976general} for details.

If we specify an initial probability distribution $P(\cdot, 0)$, the CME determines all future probability distributions $P(\cdot, t)$ for $t > 0$.
Often this initial distribution is assumed to be a Dirac distribution, $P(\cdot, 0) = \delta_{\hat{\mbf x}_0}$, where all of the probability is concentrated on a single state $\hat{\mbf x}_0 \in \hat{\mathcal X}$. 
However, in principle, the initial distribution could be supported on any subset of $\hat{\mathcal X}$.

Note that the CME holds for \emph{all} reachable states $\hat{\mbf x} \in \hat{\mathcal X}$.
So it is not just a single equation but a whole system of equations.
This system can be written concisely as 
\begin{equation} \label{ODE CME}
\frac{d\mbf p}{dt}(t) = \mbf G \mbf p(t),
\end{equation}
where $\mbf G$ is a time-invariant (infinitesimal generator) matrix whose coefficients are linked to the propensity functions, and $\mbf p$ is a vector of probabilities with one component for each $\hat {\mbf x} \in \hat{\mathcal X}$.
The initial probability distribution is now represented as $\mbf p(0)$.
While this equation is conceptually simple, there is often a huge number of reachable states $\hat{\mbf x} \in \hat{\mathcal X}$.
This means that the vector $\mbf p$ can have a very large (or even infinite) dimension, with $\mbf G$ being correspondingly large.
The result is that it is impractical to solve Equation \eqref{ODE CME} directly for stochastic chemical kinetic systems of any appreciable size.

\subsection{Moments in Stochastic Chemical Kinetics} \label{moments in SCK}
The probability distribution $P(\cdot ,t)$ can be characterized by its moments.  
In particular, for any multi-index $\mbf j=(j_1, \dots ,j_{\hat N} ) \in \mathbb N^{\hat N}$ we have a moment $\mu_{\mbf j} (t)$ defined as
\begin{equation} \label{definition of moments}
\mu_\mbf j (t) \equiv \sum_{\hat{\mbf x} \in \hat{\mathcal X}}  \hat{\mbf x}^\mbf j P(\hat{\mbf x},t) ,
\end{equation}
where the sum is over the set $\hat{\mathcal X}$ of all reachable states, and $\hat{\mbf x}^\mbf j = \prod_{k=1}^{\hat N} \hat{x}_k^{j_k}$  is a monomial.
The \emph{order} of the moment $\mu_\mbf j$ is defined as the sum $|\mathbf j| \equiv \sum_{k=1}^{\hat N} j_k$.
Notice that  the zeroth-order moment $\mu_\mbf 0(t)$ indexed by $\mbf 0 = (0,\dots,0)$ 
is simply the sum of probabilities across all reachable states,
so that $\mu_\mbf 0 (t)=1$ for all times $t$.  

A nice feature of moments is that, using just the low-order moments, we can express several quantities of interest that effectively summarize the distribution $P(\cdot, t)$.
For example, 
the first-order moment $\mu_{\mbf e_i}(t)$ indexed by $\mbf e_i=(0,\dots ,1, \dots, 0)$ is the mean molecular count for independent species $i \in \{1, \dots, \hat N \}$ at time $t$:  
\begin{equation} \label{first moment is mean of independent species}
\mu_{\mbf e_i} (t) \equiv \sum_{\hat{\mbf x} \in \hat{\mathcal X}} \hat{\mbf x}^{\mbf e_i}P(\hat{\mbf x},t)  
=
\sum_{\hat{\mbf x} \in \hat{\mathcal X}} \hat{x}_i P(\hat{\mbf x},t)  
=
\langle \hat X_i(t) \rangle
.
\end{equation}
The first-order moments can also be used with Equation \eqref{expanded expression of invariants} to express the mean molecular count for each dependent species $k \in \{1, \dots, L\}$.
In particular, if we let $\beta_{k,j}$ denote the element in the $k$th row and $j$th column of the matrix $\tilde{\mbf B} ^{-1}\hat{\mbf B}$, and $\alpha_k$ equal the $k$th component of the vector $\tilde{\mbf B}^{-1} \mbf f$, then we have
\begin{equation} \label{first order moments and mean of dependent species}
\begin{aligned}
\langle \tilde X_k(t) \rangle
&=
\sum_{\hat{\mbf x} \in \hat{\mathcal X}} \tilde x_k
P(\hat{\mbf x}, t), \\
&=
\sum_{\hat{\mbf x} \in \hat{\mathcal X}} \mbf e_k\trans \left (  
\tilde{\mbf B}^{-1} \mbf f - \tilde{\mbf B}^{-1} \hat{\mbf B} \hat{\mbf x}
\right )
P(\hat{\mbf x}, t), \\
&=
\sum_{\hat{\mbf x} \in \hat{\mathcal X}} \left (  
\alpha_k - \sum_{j=1}^{\hat N} \beta_{k,j} \hat{x}_j
\right )
P(\hat{\mbf x}, t), \\
&=
\alpha_k -  \sum_{j=1}^{\hat N} \beta_{k,j} \mu_{\mbf e_j}(t).
\end{aligned}
\end{equation}

Coming to the second-order moments, we see that $\mu_{2 \mbf e_i}(t)$ is equal to $\langle \hat X_i^2(t) \rangle$.
So, $\mu_{\mbf e_i}(t)$ and $\mu_{2 \mbf e_i}(t)$ can be used together to compute the variance in the count of molecules of independent species $i$ at time $t$: 
\begin{equation} \label{variance in terms of moments}
\sigma_{i}^2(t) \equiv \langle \hat X_i^2(t)\rangle - \langle \hat X_i(t)\rangle^2 = \mu_{2 \mbf e_i}(t) - \mu_{\mbf e_i}^2(t).
\end{equation}
Similarly, the moments can be used to compute covariances between independent species $i$ and $j$:
\begin{equation} \label{covariances between independent species}
\begin{aligned}
\sigma_{i,j}^2(t) &\equiv \langle \hat X_i(t) \hat X_j(t) \rangle - \langle \hat X_i(t)\rangle \langle \hat X_j(t)\rangle \\
&= \mu_{\mbf e_i + \mbf e_j}(t) - \mu_{\mbf e_i}(t) \mu_{\mbf e_j}(t).
\end{aligned}
\end{equation}

The appeal of working with moments is that they allow us to bypass the problem of high dimensionality that we encountered in Equation \eqref{ODE CME}.
We give up a complete description of the probability distribution $P(\cdot, t)$ in the terms of the high-dimensional vector $\mbf p(t)$ in favor of a summary description in terms of its low-order moments.
In principle, this trade-off allows us to compute properties of stochastic chemical kinetic systems for which solving the CME more directly is computationally intractable.

\subsection{The Closure Problem} \label{closure problem}
As described by Smadbeck and Kaznessis\cite{smadbeck2013closure}, Sotiropoulos and Kaznessis\cite{sotiropoulos2011analytical}, and C. S. Gillespie\cite{gillespie2009moment}, the CME can be used to derive a system of linear ordinary differential equations describing how the moments of the distribution $P(\cdot, t)$ change over time. 
For reaction systems containing at most first-order (i.e., unimolecular) reactions, things work out nicely:
we can pick an arbitrary $m \in \mathbb N$, and construct the ODE describing how the moments up to order $m$ change over time:
\begin{equation}
\frac{d\bm \mu_L}{dt}(t)=\mbf A_L\bm \mu_L(t),
\end{equation}
where $\bm \mu_L(t)$ is a vector of \lq\lq{}low-order\rq\rq{} moments order up to order $m$, and $\mbf A_L$ is a constant matrix.
However, if the reaction system contains any reactions of order $q > 1$ (e.g., bimolecular reactions), then the ODE becomes
\begin{equation} \label{closure problem ode}
\frac{d\bm \mu_L}{dt}(t)=\mbf A_L\bm \mu_L(t)+ \mbf A_H\bm \mu_H(t),
\end{equation}
where $\bm \mu_H(t)$ is a vector of \lq\lq{}high-order\rq\rq{} moments, order $m + 1$ to $m + q - 1 \equiv M$.
So the time derivatives of the low-order moments depend on high-order moments.
This is the infamous \lq\lq{}closure problem\rq\rq{}.
It is unclear how to solve such a dynamic system.

The closure problem also frustrates even a relatively simple steady-state analysis.
What we\rq{}d like to do is set the left-hand side of Equation \eqref{closure problem ode} equal to zero 
\begin{equation} \label{closure problem steady-state}
\mbf 0=\mbf A_L\bm \mu_{L,\text{ss}}+ \mbf A_H\bm \mu_{H,\text{ss}},
\end{equation}
and solve for the steady-state moments $\bm \mu_{L,\text{ss}}$  and $\bm \mu_{H,\text{ss}}$ of the steady-state probability distribution $P_\text{ss}(\cdot) \equiv \lim_{t \rightarrow + \infty} P(\cdot, t)$, assuming some specified initial distribution $P(\cdot, 0)$. 
Assuming we could calculate the vector $\bm \mu_{L,\text{ss}}$, we could extract the steady-state values of $\langle \hat X_i \rangle_\text{ss} \equiv \mu_{{\mbf e_i},\text{ss}}$ and $(\sigma_{i}^2)_\text{ss} \equiv \mu_{2 \mbf e_i,\text{ss}}- \mu_{\mbf e_i,\text{ss}}^2$ for each independent species $i$.
The trouble is that Equation \eqref{closure problem steady-state} is under-determined: it has more unknowns than linearly independent equations.
Even if we leverage our a priori knowledge of probability distributions and set $\mu_{\mbf 0,\text{ss}} = 1$, one can show there are still more unknowns than linearly independent equations.
This means that the system has infinitely many solutions, and we can\rq{}t simply solve for the steady-state moments $\bm \mu_{L,\text{ss}}$ and $\bm \mu_{H,\text{ss}}$.

\subsection{Bounds on Steady-State Systems} \label{steady state section}
In our previous paper\cite{dowdy2018bounds}, we described a paradigm for calculating bounds on quantities of interest for steady-state probability distributions.
This paradigm consisted of writing down several mathematical conditions that the steady-state moment vector $\bm \mu_\text{ss}$ must necessarily satisfy, and then optimizing over all vectors $\tilde{\bm \mu}_\text{ss}$ that satisfy these conditions, searching for that vector which maximizes or minimizes the quantity of interest.
For example, the optimization problem for calculating an upper bound on the mean molecular count of species $i$ at steady state, $\langle \hat X_i \rangle_\text{ss}$, can be written abstractly as
\begin{equation} \label{abstract optimization problem, original}
\begin{aligned}
\langle \hat X_i \rangle_\text{ss}^U \equiv
& \max_{\tilde{\bm \mu}_\text{ss}}  && \tilde \mu_{\mbf e_i,\text{ss}} \\
&\ \text{ s.t. } && \tilde{\bm \mu}_\text{ss} \text{ satisfies necessary steady-state} \\
&	&&\text{moment conditions}. 
\end{aligned}
\end{equation}
Note that we are making a distinction between the actual steady-state moment vector $\bm \mu_\text{ss}$ and the decision variables $\tilde{\bm \mu}_\text{ss}$ which serve as a proxy for $\bm \mu_\text{ss}$.

The optimal value of Problem \eqref{abstract optimization problem, original} is guaranteed to be an upper bound on the true $\langle \hat X_i \rangle_\text{ss}$, because the true steady-state moment vector $\bm \mu_\text{ss}$ is  a feasible point for the optimization problem by construction.
This reasoning is valid whether $\tilde{\bm \mu}_\text{ss}$ and $\bm \mu_\text{ss}$ are considered to be infinite sequences or vectors containing only finitely many moments.
However, for practical computations, we must work with finite vectors.
So, going forward, we will specify that $\tilde{\bm \mu}_\text{ss}$ contains only those moments up to order $2n \in \mathbb N$, where $n = \lceil \frac{M}{2}\rceil$.
The reason for this choice of $n$ is explained in our previous paper\cite{dowdy2018bounds}.

Our list of necessary conditions consisted of three main parts:
first, Equation \eqref{closure problem steady-state}, expressed in terms of the decision variables $\tilde{\bm \mu}_\text{ss}$,
\begin{equation} \label{closure problem steady-state, tildes}
\mbf 0=\mbf A_L\tilde{\bm \mu}_{L,\text{ss}}+ \mbf A_H\tilde{\bm \mu}_{H,\text{ss}};
\end{equation}
second, the fact that the total probability is one,
\begin{equation}
\tilde{\bm \mu}_{\mathbf 0, \text{ss}} = 1;
\end{equation}
and third, several linear matrix inequalities (LMIs) derived solely from the fact that the unknown probability distribution is supported on the set $\hat{\mathcal X} \subset \bar{\mathcal X}$:
\begin{equation} \label{basic moment lmi}
\mbf M_n^{\mbf 0} (\tilde{\bm \mu}_\text{ss} ) \succeq \mbf 0,
\end{equation}
\begin{equation} \label{nonegative independents}
\mbf M_{n-1}^{\mbf e_j} (\tilde{\bm \mu}_\text{ss} ) \succeq \mbf 0, \ \ \ \forall j \in \{1, \dots, \hat N\},
\end{equation}
\begin{equation} \label{nonegative dependents}
\begin{aligned}
\alpha_k \mbf M_{n-1}^{\mbf 0} (\tilde{\bm \mu}_\text{ss} )  
&-  \sum_{j=1}^{\hat N} \beta_{k,j} \mbf M_{n-1}^{\mbf e_j} (\tilde{\bm \mu}_\text{ss} ) \succeq \mbf 0, \\
&\ \ \ \forall k \in \{1, \dots, L\},
\end{aligned}
\end{equation}
The exact definitions of the matrices $\mbf M_n^{\mbf 0} (\tilde{\bm \mu}_\text{ss})$, $\mbf M_{n-1}^{\mbf 0} (\tilde{\bm \mu}_\text{ss})$, and $\mbf M_{n-1}^{\mbf e_j} (\tilde{\bm \mu}_\text{ss})$ can be found in the supplementary material of our previous publication\cite{dowdy2018bounds}.
However, the important point is that these matrices are symmetric and linear with respect to their arguments.
Each LMI simply asserts that the matrices on the left-hand side of the \lq\lq{}$\succeq$\rq\rq{} must be positive semidefinite (i.e., have all nonnegative eigenvalues).

Substituting in these necessary conditions gives us an SDP for calculating $\langle \hat X_i \rangle_\text{ss}^U$.
By changing the \lq\lq{}max\rq\rq{} to a \lq\lq{}min\rq\rq{}, we can calculate the lower bound $\langle \hat X_i \rangle_\text{ss}^U$, and by variations on this theme, we can calculate bounds on other quantites, such as the steady-state variance in the molecular count of species $i$.

Our paradigm for calculating time-varying bounds on dynamic systems will be similar.
In fact, we will make use of some of the same necessary conditions that appear above.

\section{Bounds on Dynamic Systems} \label{the bounding method section}
In this section, we extend the method for calculating bounds on the steady-state stochastic chemical kinetic systems to calculate bounds on dynamic systems.

%

\subsection{The Paradigm} \label{the paradigm}
Suppose that we have a generic stochastic chemical kinetic system, characterized by a stoichiometry matrix $\mbf S \in \mathbb Z^{N \times R}$ and a vector of rate constants
$\mbf c \in \R^R$.
Assume that there is at least one reaction with order greater than one, so that this system exhibits the closure problem when subjected to a moment analysis.
Suppose that we have analyzed $\mbf S$ to construct an invariant matrix $\mbf B \in \R^{L \times N}$, as described in Section \ref{invariants}, and that we know the associated invariant values $\mbf f \in \R^L$.
Suppose further that have identified the $\hat N = N - L$ chemical species we wish to treat as independent and constructed the matrices $\hat{\mbf B} \in \R^{L \times \hat N}$ and $\tilde{ \mbf B} \in \R^{L\times L}$.
Finally, suppose that we have chosen a value of $m \in \mathbb N$ and constructed the matrices $\mbf A_L$ and $\mbf A_H$ described in Section \ref{closure problem}.
We are interested in analyzing the properties of the probability distribution describing the stochastic chemical kinetic system at a particular time $T$.

Consider the problem of bounding $\langle \hat X_i(T) \rangle$, the mean count of molecules of independent species $i$ at time $T$.
What we\rq{}d like to do is calculate two numbers $\langle \hat X_i (T) \rangle^L$ and $\langle \hat X_i (T) \rangle^U$ such that
\begin{equation}
\langle \hat X_i (T)\rangle^L\leq \langle \hat X_i (T)\rangle \leq \langle \hat  X_i (T)\rangle^U
\end{equation}
is guaranteed.

To calculate these bounds, we will again make use of the paradigm described in Section \ref{steady state section}, only this time our necessary conditions will be on the probability distribution at time $T$, not at steady state.
In particular, the abstract problem for calculating the upper bound $\langle \hat X_i(T) \rangle^U$ is:
\begin{equation} \label{abstract optimization problem}
\begin{aligned}
\langle \hat X_i(T) \rangle^U \equiv
& \max_{\tilde{\bm \mu}(T)}  && \tilde{\mu}_{\mbf e_i}(T) \\
&\ \text{ s.t. } && \tilde{\bm \mu}(T) \text{ satisfies necessary} \\
&	&&\text{moment conditions at time } T. 
\end{aligned}
\end{equation}
As before, we make a distinction between the true moment vector $\bm \mu(T)$ at time $T$, and the decision variable $\tilde{\bm \mu}(T)$, which is a proxy for $\bm \mu(T)$.

Following the same reasoning, we can calculate a lower bound on $\langle \hat X_i (T)\rangle$ by \emph{minimizing} over the set of vectors $\tilde{\bm \mu}(T)$ satisfying the necessary moment conditions.


\subsection{Necessary Moment Conditions} \label{necessary steady-state conditions}
What exactly are the necessary moment conditions appearing in Problem \eqref{abstract optimization problem}?
As before, we must have that the total probability is equal to one:
\begin{equation} \label{sum of probabilities equals 1}
\tilde \mu_{\mbf 0}(T) = 1.
\end{equation}

Also, because the distribution $P(\cdot, T)$ is supported on the set $\hat{\mathcal X} \subset \bar{\mathcal X}$, we again have the LMIs that were relevant in the steady-state analysis:

\begin{equation} \label{basic moment lmi}
\mbf M_n^{\mbf 0} (\tilde{\bm \mu}(T)) \succeq \mbf 0,
\end{equation}
\begin{equation} \label{nonegative independents}
\mbf M_{n-1}^{\mbf e_j} (\tilde{\bm \mu}(T)) \succeq \mbf 0, \ \ \ \forall j \in \{1, \dots, \hat N\},
\end{equation}
\begin{equation} \label{nonegative dependents}
\begin{aligned}
\alpha_k \mbf M_{n-1}^{\mbf 0} (\tilde{\bm \mu}(T))  
&-  \sum_{j=1}^{\hat N} \beta_{k,j} \mbf M_{n-1}^{\mbf e_j} (\tilde{\bm \mu}(T)) \succeq \mbf 0, \\
&\ \ \ \forall k \in \{1, \dots, L\}.
\end{aligned}
\end{equation}
The set of vectors satisfying LMIs \eqref{basic moment lmi}-\eqref{nonegative dependents} is a mathematical cone.
To simplify the notation in what follows, will represent this cone concisely as $C_n(\bm \alpha, \bm \beta)$.
Thus, 
\begin{equation} \label{mu is in cone}
\tilde{\bm \mu}(T) \in C_n(\bm \alpha, \bm \beta)
\end{equation}
is equivalent to LMIs \eqref{basic moment lmi}-\eqref{nonegative dependents}.

Conditions \eqref{sum of probabilities equals 1}-\eqref{nonegative dependents} are notably lacking any information about the dynamics of the system.
To obtain necessary conditions implied by the dynamics, we make use of Equation \eqref{closure problem ode}, which holds for all times $t$.
Suppose that we pick an arbitrary $\rho \in \mathbb R$, multiply both sides of Equation \eqref{closure problem ode} by $e^{\rho(T - t)}$, and then integrate from $t = 0$ to $t = T$:
\begin{equation} \label{closure problem ode integrated}
\begin{aligned}
\int_{0}^T &e^{\rho(T - t)} \frac{d\bm \mu_L}{dt}(t)dt \\
&=
\int_{0}^T e^{\rho(T - t)} (\mbf A_L\bm \mu_L(t)+ \mbf A_H\bm \mu_H(t) )dt.
\end{aligned}
\end{equation}
Applying integration by parts to the left-hand side, we obtain
\begin{equation} \label{closure problem ode integrated, LHS}
\begin{aligned}
\int_{0}^T &e^{\rho(T - t)} \frac{d\bm \mu_L}{dt}(t)dt 
\\
&= e^{\rho(T - t)} \bm \mu_L(t) |_{0}^T - \int_0^T (-\rho) e^{\rho(T - t)} \bm \mu_L(t) dt,
\\
&= \bm \mu_L(T) -  e^{\rho T} \bm \mu_L(0) + \rho \int_0^T e^{\rho(T - t)} \bm \mu_L(t) dt.
\end{aligned}
\end{equation}
We presume that the initial values of the low-order moments $\bm \mu_L(0)$ can be easily computed from the initial distribution $P(\cdot, 0)$ via Equation \eqref{definition of moments}.
This is true, for example, if the initial molecular count is known exactly -- which corresponds to an initial probability distribution $P(\cdot, 0)$ where all the probability is concentrated on a single state $\hat{\mbf x}_0$, i.e., the Dirac distribution $\delta_{\hat{\mbf x}_0}$.
However, it may also be the case that we \emph{don\rq{}t} know the initial molecular count exactly.
In this case, our initial probability distribution $P(\cdot, 0)$ will be supported on several reachable states $\hat{\mbf x} \in \hat{\mathcal X}$.
Our method can handle this situation, as long as we can compute the moments $\bm \mu_L(0)$ (see Section \ref{uncertain initial state section}).

For the right-hand side, we can make use of the fact that the integral is a linear operator to obtain
\begin{equation} \label{closure problem ode integrated, RHS}
\mbf A_L \int_{0}^T e^{\rho(T - t)}\bm \mu_L(t) dt + \mbf A_H \int_{0}^T e^{\rho(T - t)}\bm \mu_H(t) dt.
\end{equation}
If we define 
\begin{equation} \label{definition of z variables}
\begin{aligned}
\mathbf z_L^{(\rho)} &\equiv \int_{0}^T e^{\rho(T - t)}\bm \mu_L(t) dt, \\
\mathbf z_H^{(\rho)} &\equiv \int_{0}^T e^{\rho(T - t)}\bm \mu_H(t) dt,
\end{aligned}
\end{equation}
we can express Equation \eqref{closure problem ode integrated} concisely as
\begin{equation}
 \bm \mu_L(T) - e^{\rho T}\bm \mu_L(0) + \rho \mathbf z_L^{(\rho)}
=
\mbf A_L \mathbf z_L^{(\rho)} + \mbf A_H \mathbf z_H^{(\rho)}.
\end{equation}
Rearranging, we obtain
\begin{equation} \label{dynamic equation for moments}
\bm \mu_L(T) - e^{\rho T} \bm \mu_L(0)
=
(\mbf A_L - \rho \mathbf I) \mathbf z_L^{(\rho)} + \mbf A_H \mathbf z_H^{(\rho)}.
\end{equation}

As before, we will replace the unknown $\bm \mu_L(T)$ with its decision variable proxy $\tilde{\bm \mu}_L(T)$.
Similarly, the vectors $\mathbf z_L^{(\rho)}$ and $\mathbf z_H^{(\rho)}$ are also unknown and will be replaced with decision variable proxies  $\tilde{\mathbf z}_L^{(\rho)}$ and $\tilde{\mathbf z}_H^{(\rho)}$, respectively.  
So necessary condition \eqref{dynamic equation for moments} becomes the following constraint in our optimization problem:
\begin{equation} \label{dynamic equation for moments, tildes}
\tilde{\bm \mu}_L(T) - e^{\rho T} \bm \mu_L(0)
=
(\mbf A_L - \rho \mathbf I) \tilde{\mathbf z}_L^{(\rho)} + \mbf A_H \tilde{\mathbf z}_H^{(\rho)}.
\end{equation}

Now, by itself, Equation \eqref{dynamic equation for moments, tildes} isn\rq{}t very useful as a constraint on $\tilde{\bm \mu}(T)$, because it is in terms of the unknown vector $\tilde{\mbf z}^{(\rho)} \equiv (\tilde{\mbf z}_L^{(\rho)}, \tilde{\mbf z}_H^{(\rho)})$.
It tells us only that $\tilde{\bm \mu}_L(T) - e^{\rho T} \bm \mu_L(0)$ must be contained in the column space of the matrix $[(\mbf A_L - \rho \mathbf I)\ \  \mathbf A_H]$.
However, if we can constrain the set of possible $\tilde{\mathbf z}^{(\rho)}$ values, Equation \eqref{dynamic equation for moments, tildes} is more useful.
To do this, we return to LMIs \eqref{basic moment lmi} - \eqref{nonegative dependents}, written for the true moment vector $\bm \mu (T)$.
Since these LMIs are derived solely from the fact that the unknown probability distribution is supported on $\hat{\mathcal X} \subset \bar{\mathcal X}$, they hold not just at time $T$, but also for all times $t \in [0,T]$.
For example, we have
\begin{equation} \label{basic moment lmi, all t}
\mbf M_n^{\mbf 0} (\bm \mu(t)) \succeq \mbf 0, \ \ \ \forall t \in [0, T].
\end{equation}
Multiplying both sides of the LMI by the nonnegative factor $e^{\rho(T - t)}$ and integrating over $[0, T]$ maintains the LMI:
\begin{equation} \label{basic moment lmi, all t, integrated}
\int_0^T e^{\rho(T - t)} \mbf M_n^{\mbf 0} (\bm \mu(t)) dt \succeq \mbf 0.
\end{equation}
Furthermore, because the integral is a linear operator, and because $\mbf M_n^{\mbf 0} (\cdot)$ is a linear function of its argument, we can bring the integral inside:
\begin{equation} \label{basic moment lmi, all t, integrated, integral inside}
\mbf M_n^{\mbf 0} \left (\int_0^T e^{\rho(T - t)} \bm \mu(t)dt \right )  = \mbf M_n^{\mbf 0} \left (\mathbf z^{(\rho)}\right ) \succeq \mbf 0.
\end{equation}
Following similar reasoning, we can show that
\begin{equation} \label{nonegative independents, z}
\mbf M_{n-1}^{\mbf e_j} (\bm z^{(\rho)}) \succeq \mbf 0, \ \ \ \forall j \in \{1, \dots, \hat N\},
\end{equation}
\begin{equation} \label{nonegative dependents,z}
\begin{aligned}
\alpha_k \mbf M_{n-1}^{\mbf 0} (\bm z^{(\rho)})  
&-  \sum_{j=1}^{\hat N} \beta_{k,j} \mbf M_{n-1}^{\mbf e_j} (\bm z^{(\rho)}) \succeq \mbf 0, \\
&\ \ \ \forall k \in \{1, \dots, L\}.
\end{aligned}
\end{equation}
LMIs \eqref{basic moment lmi, all t, integrated, integral inside} - \eqref{nonegative dependents,z} can be written concisely as
\begin{equation} \label{zp is in cone too}
\bm z^{(\rho)} \in C_n(\bm \alpha, \bm \beta).
\end{equation}
We have shown that membership in the cone $C_n(\bm \alpha, \bm \beta)$ is a necessary condition for the vector $\bm z^{(\rho)}$.
Accordingly, we will enforce this membership as a constraint on its decision variable proxy $\tilde{\bm z}^{(\rho)}$:
\begin{equation} \label{zp is in cone too, tilde}
\tilde{\bm z}^{(\rho)} \in C_n(\bm \alpha, \bm \beta).
\end{equation}

Recall that our choice of $\rho \in \R$ was arbitrary.  
It follows that conditions \eqref{dynamic equation for moments} and \eqref{zp is in cone too} can be written for any $\rho \in \R$.
In fact, they hold for each $\rho$ in any subset $\mathcal R \subset \R$.
It follows that we can write constraints \eqref{dynamic equation for moments, tildes} and \eqref{zp is in cone too, tilde} for each $\rho$ in any subset $\mathcal R \subset \R$.


\subsection{A Semidefinite Program}
If we use constraints \eqref{sum of probabilities equals 1}, \eqref{mu is in cone}, \eqref{dynamic equation for moments, tildes}, and \eqref{zp is in cone too, tilde} in place of the abstract statement \lq\lq{}$\tilde{\bm \mu}(T)$ satisfies necessary moment conditions at time $T$\rq\rq{}, we obtain Optimization Problem \eqref{concrete optimization problem}:

\begin{equation} \label{concrete optimization problem}
\begin{aligned}
\langle \hat X_i(T)\rangle^U =& \max_{\substack{\tilde{\bm \mu}(T), \\ \tilde{\mathbf z}^{(\rho)}, \forall \rho \in \mathcal R}}  && \tilde \mu_{\mbf e_i}(T) \\
&\ \  \ \text{ s.t. } && \tilde \mu_{\mbf 0}(T) = 1, \\
& && \tilde{\bm \mu}(T) \in C_n(\bm \alpha, \bm \beta), \\
& && \tilde{\mbf z}^{(\rho)} \in C_n(\bm \alpha, \bm \beta), \ \ \ \forall \rho \in \mathcal R, \\
& && \text{Equation \eqref{dynamic equation for moments, tildes} holds},  \ \ \ \forall \rho \in \mathcal R.
\end{aligned}
\end{equation}

Note that the vectors $\tilde{\mbf z}^{(\rho)}$ for all $\rho \in \mathcal R$ are decision variables in addition to the vector $\tilde{\bm \mu}(T)$.
As with the vector $\tilde{\bm \mu}(T)$, it is only necessary for these vectors to contain moments up through order $2n$, where $n \equiv \lceil \frac{M}{2} \rceil$.

With its linear objective function, linear equations, and LMIs, Problem \eqref{concrete optimization problem} is a special type of optimization problem called a Semidefinite Program (SDP).
As with all SDPs, Problem \eqref{concrete optimization problem} is convex.
Thus, at least in theory, we should be able to solve it efficiently\cite{vandenberghe1996semidefinite}.
Doing so, we obtain the desired upper bound, $\langle \hat X_i(T)\rangle^U$.
Solving the corresponding minimization problem, we obtain the lower bound, $\langle \hat X_i(T)\rangle^L$.

\subsection{Inspiration from Previous Work} \label{previous work section}
The inspiration for the bounding method described in the preceding sections comes from a paper by Bertsimas and Caramanis\cite{bertsimas2006bounds}, in which moment-based SDPs are used to bound the solutions of linear partial differential equations (PDEs).
The central idea of their method is to view the solution $u(\cdot)$ of the PDE as a distribution over the problem domain $\Omega$.
Taking this view, they define the \emph{full moments}
\begin{equation} \label{full moments definition}
m_\mbf j \equiv \int_\Omega \mbf x^{\mbf j} u(\mbf x)
\end{equation}
and  \emph{boundary moments}
\begin{equation}
b_\mbf j \equiv \int_{\partial \Omega} \mbf x^{\mbf j} u(\mbf x)
\end{equation}
of the distribution, where $\partial \Omega$ is some portion of the boundary.
Starting from the linear PDE and the associated boundary conditions, they derive linear equations that these moments must satisfy.
Furthermore, they derive LMIs that the moments must satisfy, simply by virtue of being moments of a distribution supported on $\Omega$.
They then solve an SDP to optimize over all vectors ($\mbf m$, $\mbf b$) which satisfy these necessary conditions, searching for that vector which maximizes or minimizes some moment of interest.

Clearly, this is thematically similar to the bounding method we have proposed for stochastic chemical kinetic systems. 
We now elaborate on this connection.
In considering the problem of stochastic chemical kinetics, we naturally focus on $P(\cdot, t)$ as a probability distribution over the reachable states $\hat{\mathcal X}$ for each time $t \in [0,T]$.
However, we can also think of the function $P(\cdot, \cdot)$ as a generalized distribution over both state space and time -- that is, a distribution supported on the set $\Omega = \hat{ \mathcal X} \times [0,T]$.
This $P(\cdot, \cdot)$ is directly analogous to the function $u(\cdot)$ above.
Furthermore, the moments $\mu_{\mbf j}(0)$ and $\mu_{\mbf j}(T)$ are analogous to the \lq\lq{}boundary moments\rq\rq{}, as they are   associated with the boundaries of $\Omega$ corresponding to $t = 0$ and $t = T$.
Finally, the quantities $z_{\mbf j}^{(\rho)}$ are analogous to the \lq\lq{}full moments\rq\rq{} above.


This last analogy may not be so obvious, but it becomes clearer if we expand Equation \eqref{definition of z variables} using Equation \eqref{definition of moments}.
Doing so, we see that for any $\mbf j \in \mathbb N^{\hat N}$,
\begin{equation}
z_{\mbf j}^{(\rho)} = \int_0^T \sum_{\hat{\mbf x} \in \hat{\mathcal X}} e^{\rho(T-t)} \hat{\mbf x}^{\mbf j} P(\hat{\mbf x}, t) dt,
\end{equation}
which can be written more abstractly, closer to Bertsimas and Caramanis\rq{}s notation, as
\begin{equation} \label{definition of z}
z_{\mbf j}^{(\rho)} = \int_\Omega e^{\rho(T-t)} \hat{\mbf x}^{\mbf j} P(\hat{\mbf x}, t),
\end{equation}
where, again, $\Omega = \hat{\mathcal X} \times [0,T]$.
When the equation for $z_{\mbf j}^{(\rho)}$ is written in this form, the analogy with Equation \eqref{full moments definition} is obvious.

The reader might protest that a closer analogy to Equation \eqref{full moments definition} would be
\begin{equation} \label{closer z analogy}
z_{\mbf j}^{(\rho)} = \int_\Omega t^\rho \hat{\mbf x}^{\mbf j} P(\hat{\mbf x}, t),
\end{equation}
and we agree.
Our departure from the strict analogy is deliberate.
As Bertsimas and Caramanis point out, while moments are classically defined in terms of monomials, we are free to define them in terms of other basis functions which may be better suited to the problem at hand.
This is exactly what we have done in our definition of  $z_{\mbf j}^{(\rho)}$.
Recall that the CME \eqref{ODE CME} is a linear time-invariant ODE:
\begin{equation*} \label{ODE CME, restated}
\frac{d\mbf p}{dt}(t) = \mbf G \mbf p(t).
\end{equation*}
Assuming that the number of reachable states $|\hat{\mathcal X}|$ is finite, and assuming that $\mbf G$ has $|\hat{\mathcal X}|$ distinct eigenvalues $\{\lambda_j\}_{j = 1}^{|\hat{\mathcal X}|}$, the solution to this system can be written as
\begin{equation}
\mbf p(t) = \sum_{j = 1}^{|\hat{\mathcal X}|} a_j e^{\lambda_j t} \mbf v_j,
\end{equation}
where the $\{\mbf v_j\}_{j = 1}^{|\hat{\mathcal X}|}$ are the right eigenvectors of $\mbf G$, and the $\{a_j\}_{j = 1}^{|\hat{\mathcal X}|}$ are complex-valued coefficients derived from the initial distribution $\mbf p(0)$.
In this case, the solution\rq{}s time-variation has an exponential character.
This strongly suggests that, in our efforts to bound the solution, we should use basis functions which are also exponential with respect to time.
Furthermore, it strongly suggests that the coefficients $\rho$ appearing in these basis functions should match the eigenvalues of the matrix $\mbf G$.






\subsection{Choosing the Values of $\rho$} \label{sec: choosing rho values}
An obvious problem with the idea of choosing our values of $\rho$ to match the eigenvalues of $\mbf G$ is that there can be as many distinct eigenvalues as there are reachable states -- often a huge number.
Recall that each value of $\rho \in \mathcal R$ has an associated collection of decision variables $\tilde{\mbf z}^{(\rho)}$ and constraints in SDP \eqref{concrete optimization problem}.
It is not tractable to have such a large number of variables and constraints; so we can only hope to use some relatively small subset of the eigenvalues in defining the set $\mathcal R$.

This brings us to the question: which eigenvalues should we use?
Our computational experience suggests that we should pick the values of $\rho$ to approximate the real parts of the first several distinct eigenvalues of the matrix $\mbf G$ when listed in order of increasing magnitude.
By the construction of $\mbf G$, one of these eigenvalues is guaranteed to be zero,
so we will always have $\rho = 0$ as one of our members of $\mathcal R$.
Using Gershgorin\rq{}s Circle Theorem\cite{gershgorin1931uber}, one can show that the nonzero eigenvalues of $\mathbf G$ all have strictly negative real parts.


The next question is: how can we calculate the eigenvalues we\rq{}d like to use in defining the set $\mathcal R$?
Since the matrix $\mbf G$ is large and sparse, an iterative Krylov subspace method\cite{saad1992numerical}  seems appropriate.
However, the fact that $\mbf G$ can be infinitely large means that the standard algorithms cannot be applied without some modification.
We are developing a modified, infinite-dimensional Krylov method, which will be a subject of a future publication.


\subsection{Bounds on the Variance}
As explained in our previous paper\cite{dowdy2018bounds}, through some relatively simple modifications of the SDP for calculating bounds on the steady-state mean molecular count of species $i$, we can construct an SDP for calculating an upper bound on the variance in this count.
The same reasoning applies for the dynamic problem, giving us the following SDP for calculating an upper bound on the variance in the molecular count of species $i$ at time $T$:


\begin{equation} \label{concrete optimization problem, variance}
\begin{aligned}
\sigma_i^2(T)^U =& \max_{\substack{\tilde{\bm \mu}(T), s, \\ \tilde{\mathbf z}^{(\rho)}, \forall \rho \in \mathcal R}}  && s \\
&\ \  \ \text{ s.t. } &&
\left [
\begin{array}{cc}
\tilde \mu_{2\mbf e_{i}}(T) - s & \tilde \mu_{\mbf e_{i}} \\
\tilde \mu_{\mbf e_{i}} & 1
\end{array}
\right ] \succeq \mbf 0,\\
& && \tilde \mu_{\mbf 0}(T) = 1, \\
& && \tilde{\bm \mu}(T) \in C_n(\bm \alpha, \bm \beta), \\
& && \tilde{\mbf z}^{(\rho)} \in C_n(\bm \alpha, \bm \beta), \ \ \ \forall \rho \in \mathcal R, \\
& && \text{Equation \eqref{dynamic equation for moments, tildes} holds},  \ \ \ \forall \rho \in \mathcal R.
\end{aligned}
\end{equation}

\subsection{Bounds on Probability}
In our previous paper\cite{dowdy2018bounds}, we also formulated SDPs for calculating an upper bound on the steady-state probability that the molecular count of species $i$ is an arbitrary interval $[x_\text{min}, x_\text{max}]$, and we saw that this led to bounding histograms.
We also noted that we could bound the probability that the steady-state probability distribution assigns to an arbitrary basic semi-algebraic set, i.e., a set of the form
\begin{equation}
\{ \mbf{\hat x} \in \R^{\hat N} : g_j(\mbf{\hat x}) \geq 0, \ j = 1, \dots, K \},
\end{equation}
where each $g_j(\cdot)$ for $j = 1, \dots, K$ is a polynomial is $\mathbf{\hat{x}}$.

While we do not discuss the details here, these ideas could be extended to the dynamic problem.
For example, we could calculate an upper bound on the histogram describing the unknown probability distribution at time $T$.
Furthermore, we could bound the probability that this distribution assigns to an arbitrary basic semi-algebraic set.

%

\subsection{Conservatism in the Bounds} \label{conservatism section}
As described in our previous paper\cite{dowdy2018bounds}, there are several sources of conservatism in the bounds calculated by solving SDP \eqref{concrete optimization problem} (and its variations).
The first of these is related to the fact that our choice of $m$, the cut-off of what we consider to be a \lq\lq{}low-order\rq\rq{} moment, is somewhat arbitrary.
The second source of conservatism is that the necessary conditions appearing in SDP \eqref{concrete optimization problem} in no way reflect the physical constraint that the number of molecules of each species must be an integer.
These sources of conservatism are discussed at length in our previous paper, and the interested reader is referred there for further details.

There is, however, one source of conservatism which cannot be found in our previous paper and which is unique to the dynamic problem.
This conservatism comes from our choice of the set $\mathcal R$.
As we\rq{}ve already pointed out, Conditions \eqref{dynamic equation for moments} and \eqref{zp is in cone too} hold for all $\rho \in \R$.
However, for Problem \eqref{concrete optimization problem} to be computationally tractable, we can only enforce these conditions for some finite subset $\mathcal R \subset \R$.
In a sense, we are thus relaxing Conditions \eqref{dynamic equation for moments} and \eqref{zp is in cone too} for all $\rho \in \R$ such that $\rho \notin \mathcal R$.
Doing so may introduce some conservatism in the resulting bounds.
This suggests that adding elements $\rho \in \R$ to our set $\mathcal R$ will improve the quality of the bounds.
When we come to the examples in Section \ref{mean bounding examples}, we will see that this is, in fact, the case.

\subsection{Scaling} \label{scaling section}

As pointed out in our previous paper\cite{dowdy2018bounds}, one shortcoming of moment-based SDPs such as Problem \eqref{concrete optimization problem} is that they can give solvers numerical difficulties.
This is especially true if the SDPs are not appropriately scaled.
We discuss some strategies for scaling in our previous paper\cite{dowdy2018bounds}, so we will not go into details here.
However, we do wish to point out that, if one solves a sequence of bounding problems for increasing times $T_j$, the bounds at time $T_j$ could be helpful in appropriately scaling the problem for time $T_{j + 1}$.

\section{Toy Example} \label{mean bounding examples}
%

In this section, we apply SDPs \eqref{concrete optimization problem} and \eqref{concrete optimization problem, variance} to a simple stochastic chemical kinetic systems as a proof of concept.
Consider the simple irreversible reaction
\begin{equation} \label{Ex1 system}
\text{A} + \text{B} \autorightarrow{$c_1$}{} \text{C}
\end{equation}
with rate constant $c_1 = 1 \ \text{s}^{-1}$, and known initial molecular counts of A $=3$, B $=4$, and C $=0$.
If we select A as the species to consider independent, this translates to an initial probability distribution $P(\cdot, 0) =\delta_3$, where all of the probability is concentrated on the reduced state $\hat x = 3$.
Given that this system features a bimolecular reaction, it exhibits the closure problem when subjected to a moment analysis.

\subsection{Mean and Variance Bounds}

If we repeatedly solve SDP \eqref{concrete optimization problem} and its minimization counterpart for this system, taking $\mathcal R = \{0, -2\}$ and $m = 3$, we obtain time-varying bounds on the mean molecular counts of each species.
Similarly, if we repeatedly solve SDP \eqref{concrete optimization problem, variance} for this system, with the same $\mathcal R$ and $m$, we can obtain time-varying upper bounds on variance for each molecular count.
These bounds are shown in the top and bottom panels, respectively, of Figure \ref{fig: Ex1 Mean and Variance Bounds, reduced P}.
For comparison, we have also included the analytical means and variances provided by McQuarrie \cite{mcquarrie1967stochastic}. 

As expected, the mean bounds do, indeed, enclose the analytical means; and the variance upper bound does, indeed, exceed the analytical variance for all times $t$.
This is consistent with the theory of Section \eqref{the bounding method section}.


\begin{figure} 
    \includegraphics[scale = 0.7]{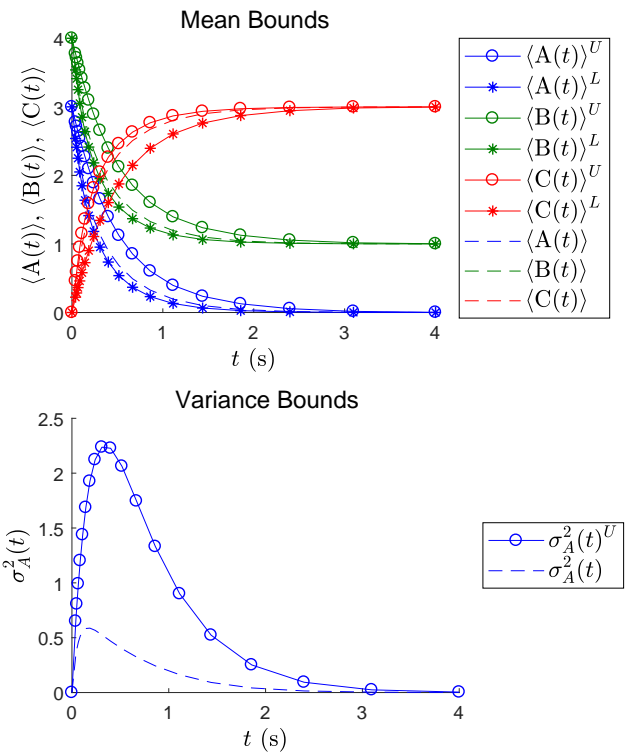}
    \caption{\label{fig: Ex1 Mean and Variance Bounds, reduced P} Time-varying bounds on System \eqref{Ex1 system}, calculated using $\mathcal R = \{0, -2\}$. The points marked with circles and stars each correspond to the solution of an SDP and are theoretically guaranteed bounds.  The lines interpolated between these points are not guaranteed bounds.  They are included just to lead the eye. The dashed lines are the analytical values, provided by McQuarrie \cite{mcquarrie1967stochastic}. The top plot shows bounds on the mean molecular count of each species.  The bottom plot shows an upper bound on the variance in the molecular count of species A.  The other species are omitted, because their variances are identical.}
\end{figure}

\subsection{Using more values of $\rho$}
In Section \ref{conservatism section}, we noted that the choice of the set $\mathcal R$ can affect the quality of the resulting bounds.
We demonstrate this by recalculating the bounds shown in Figure \ref{fig: Ex1 Mean and Variance Bounds, reduced P} with the enlarged set $\mathcal R = \{ 0, -2, -6 \}$.
The results are shown in Figure \ref{fig: Ex1 Mean and Variance Bounds, more P}.
The bounds are noticeably tighter for both the means and the variance, which is consistent with our prior reasoning.

\begin{figure} 
    \includegraphics[scale = 0.7]{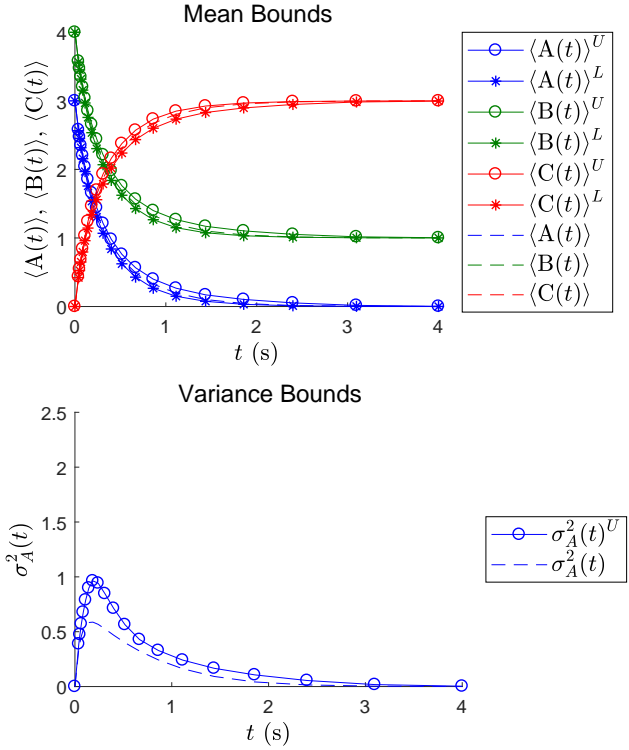}
    \caption{\label{fig: Ex1 Mean and Variance Bounds, more P} This figure is equivalent to Figure \ref{fig: Ex1 Mean and Variance Bounds, reduced P} in every way, except that the bounds were calculated using the enlarged set $\mathcal R = \{ 0, -2, -6 \}$, giving tighter bounds.}
\end{figure}

\section{A Bit More Complexity} \label{mean bounding examples}
%

In this section, we apply SDPs \eqref{concrete optimization problem} and \eqref{concrete optimization problem, variance} to a slightly more complex reaction system, where we\rq{}ve added a reversible reaction:
\begin{equation} \label{Ex6 system}
\text{A} + \text{B} \autorightarrow{$c_1$}{} \text{C} \autorightleftharpoons{$c_2$}{$c_3$} \text{D}
\end{equation}
The rate constants for this system are $c_1 = 1 \ \text{s}^{-1}, c_2 = 2.1 \ \text{s}^{-1},$ and $c_3 = 0.3 \ \text{s}^{-1}$
The initial molecular counts of A $=3$, B $=4$, C $=0$, and D $=0$.
As with the previous example, this reaction system exhibits the closure problem when subjected to a moment analysis.

\subsection{Mean and Variance Bounds}

If we repeatedly solve SDP \eqref{concrete optimization problem} and its minimization counterpart for this system, taking $\mathcal R = \{0, -2, -2.4\}$ and $m = 3$, we obtain time-varying bounds on the mean molecular counts of each species.
Similarly, if we repeatedly solve SDP \eqref{concrete optimization problem, variance} for this system, with the same $\mathcal R$ and $m$, we can obtain time-varying upper bounds on variance for each molecular count.
These bounds are shown in the top and bottom panels, respectively, of Figure \ref{fig: Ex6 Mean and Variance Bounds, reduced P}.

In this case, we have no analytical solution.
However, the plotted curves match what we would expect.
The count of molecules of species B decreases to 1, at which point the molecules of A (not shown in the plot) are exhausted.
This is the same behavior we saw in Figures \ref{fig: Ex1 Mean and Variance Bounds, reduced P} and \ref{fig: Ex1 Mean and Variance Bounds, more P}, and this makes sense, because the addition of the reversible reaction in System \eqref{Ex6 system} does not change the dynamics of species A and B.
The mean molecular count for species C rises and then falls, leveling off at about 0.5 molecules, while the mean molecular count of species D increases monotonically, leveling off at about 2.5 molecules.
The upper bounds on the variances for the two species both approach the same limiting value (about 0.75).
This makes sense, because, as the system comes to equilibrium, when no molecules of A and B remain, the probability will be distributed entirely between species C and D, and the uncertainty in the molecular count of one is equal to the uncertainty in the molecular count of the other.


\begin{figure} 
    \includegraphics[scale = 0.7]{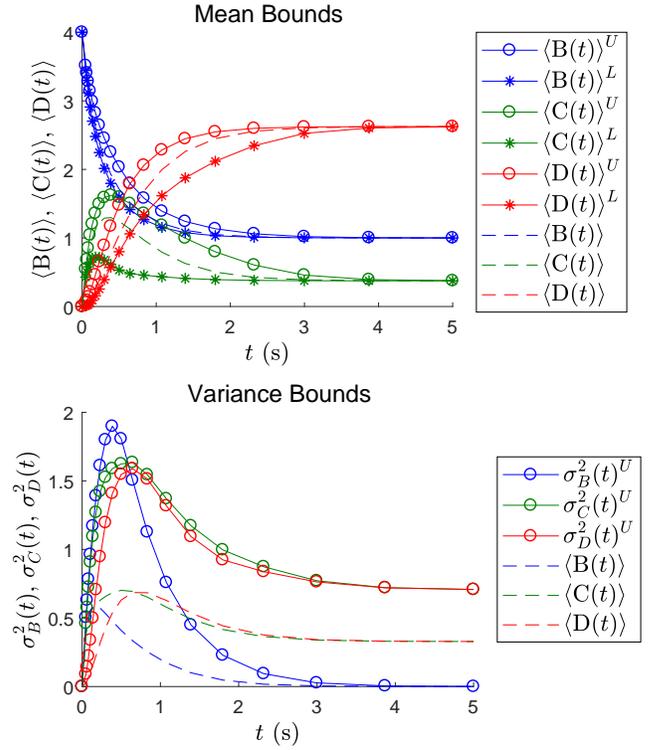}
    \caption{\label{fig: Ex6 Mean and Variance Bounds, reduced P} Time-varying bounds on System \eqref{Ex6 system}, calculated using $\mathcal R = \{0, -2, -2.4\}$. The top plot shows bounds on mean molecular counts, while he bottom plot shows upper bounds on the variances in the counts.  Species A is omitted because its behavior closely follows that of species B, and adding extra curves would only clutter the plot.  The dashed lines are the true mean and variance trajectories, obtained by direct solution of the CME.}
\end{figure}

\subsection{Using more values of $\rho$}
As with the previous example, we now add a value of $\rho$ to our set $\mathcal R$, repeat the bounding calculation, and see an improvement in the bounds.
In particular, adding $\rho = -4.4$ to our $\mathcal R$, we obtain the bounds shown in Figure \ref{fig: Ex6 Mean and Variance Bounds, more P}.
Comparing with Figure \ref{fig: Ex6 Mean and Variance Bounds, reduced P}, we see substantial improvement in the lower bound of the mean molecular count for species C.
We also see that the limiting value of the variance upper bound for species C and D is about half of its previous value.
Finally, for each species, the peak in the variance upper bound (around 0.5 s) has been reduced.

\begin{figure} 
    \includegraphics[scale = 0.7]{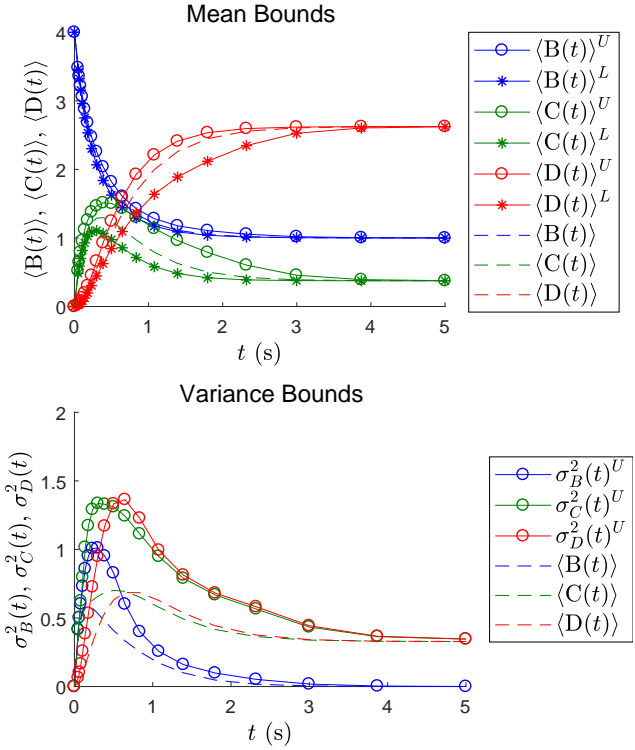}
    \caption{\label{fig: Ex6 Mean and Variance Bounds, more P} This figure is equivalent to Figure \ref{fig: Ex6 Mean and Variance Bounds, reduced P} in every way, except that the bounds were calculated using the enlarged set $\mathcal R = \{ 0, -2, -2.4, -4.4 \}$, giving better results.}
\end{figure}

\subsection{Sensitivity of the Values of $\rho$}
As explained in Section \ref{sec: choosing rho values}, while any values of $\rho$ will result in theoretically-guaranteed bounds, we recommend picking the values of $\rho$ to match the real parts of the first several distinct eigenvalues of the matrix $\mbf G$, when these eigenvalues are listed in order of increasing magnitude.
This is exactly how we chose the values of $\rho$ for the two foregoing examples.
These two examples are small enough that we can calculate the eigenvalues directly.
However, this will not be the case in general.
Usually, the best we can hope for is some numerical approximation of the eigenvalues.
This begs the question: how robust is our bounding method to the choice of $\rho$ values?
If the values of $\rho$ are off by a little bit, do the bounds become so conservative that they are practically useless?

To explore this idea, we repeated the bounding calculation for Reaction System \eqref{Ex6 system}, using a set of perturbed $\rho$ values: $\mathcal R = \{ 0, -1.9, -2.6, -4.7 \}$.
The resulting bounds are shown in Figure \ref{fig: Ex6 Mean and Variance Bounds, perturbed P}.
Comparing this plot with Figure \ref{fig: Ex6 Mean and Variance Bounds, more P}, we see that the perturbation of the values of $\rho$ did not substantially affect the quality of the computed bounds.
We see that the perturbed $\rho$ values create a slight long-time gap in the mean bounds for species C and D, which is undesirable.
However, mean bounds on these species at intermediate times (e.g., $t = 1$ s) actually seem a little tighter.
This demonstrates that the bounding method does not require exact knowledge of the eigenvalues of the underlying CME to obtain reasonable results.

That being said, the choice of $\rho$ values does matter.
Using a set of further perturbed $\rho$ values ($\mathcal R = \{ 0, -6, -12, -18 \}$), we produced the bounds shown in Figure \ref{fig: Ex6 Mean and Variance Bounds, stupid P}.
In this Figure, we see wide gaps in the long-time mean bounds for all species.
Furthermore, the variance bounds are much less tight.

In summary, while the chosen values of $\rho$ do not have to match low-magnitude eigenvalues $\mbf G$ exactly, at least approximating them seems to be a good heuristic.



%

\begin{figure} 
    \includegraphics[scale = 0.7]{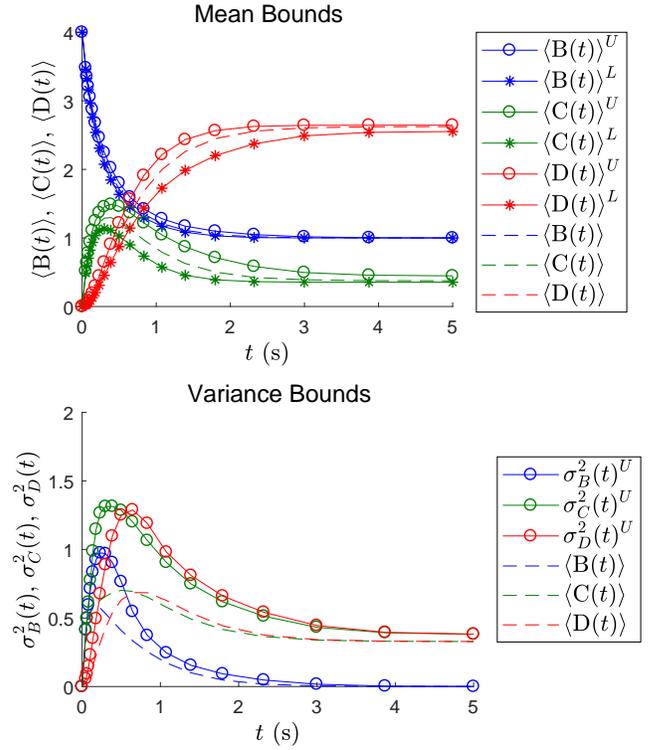}
    \caption{\label{fig: Ex6 Mean and Variance Bounds, perturbed P} This figure is equivalent to Figure \ref{fig: Ex6 Mean and Variance Bounds, more P}, except that the bounds were calculated using the perturbed set $\mathcal R = \{ 0, -1.9, -2.6, -4.7 \}$, giving slightly different results.  In particular, notice the long-time gap that has appeared in the mean bounds for species C and D.}
\end{figure}


\begin{figure} 
    \includegraphics[scale = 0.7]{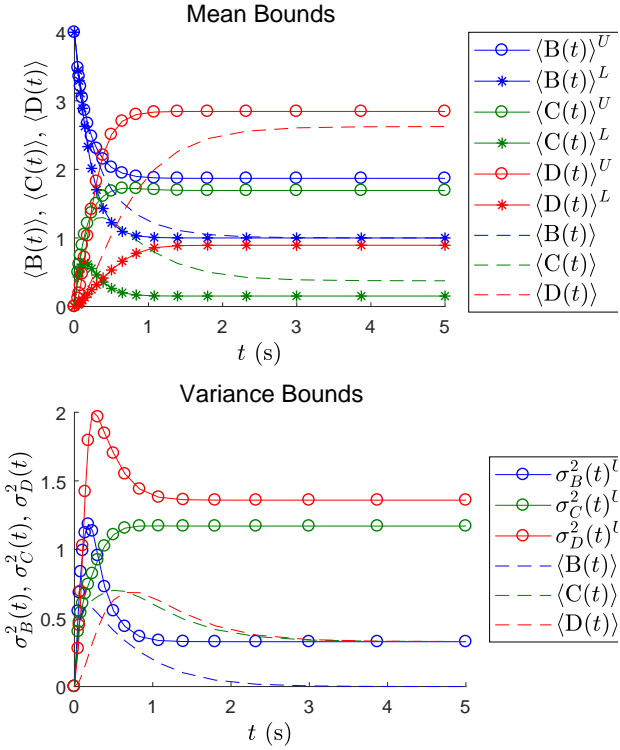}
    \caption{\label{fig: Ex6 Mean and Variance Bounds, stupid P} This figure is equivalent to Figures \ref{fig: Ex6 Mean and Variance Bounds, more P} and \ref{fig: Ex6 Mean and Variance Bounds, perturbed P}, except that the bounds were calculated using the further perturbed set $\mathcal R = \{ 0, -6, -12, -18 \}$.  This perturbation dramatically degrades the quality of the bounds.}
\end{figure}


\section{Complex Eigenvalues} \label{complex eigenvalue section}
Given our observation in Section \ref{previous work section} that it seems reasonable to choose the values of $\rho$ to match the eigenvalues of the matrix $\mbf G$, it may seem odd that, in Section \ref{sec: choosing rho values}, we suggested focusing on only the real parts of these eigenvalues.
In fact, if we know that some of the low-magnitude eigenvalues have nonzero imaginary parts, we can use this information to obtain tighter bounds.

\begin{figure}[H] 
    \centering
    \includegraphics[scale = 0.75]{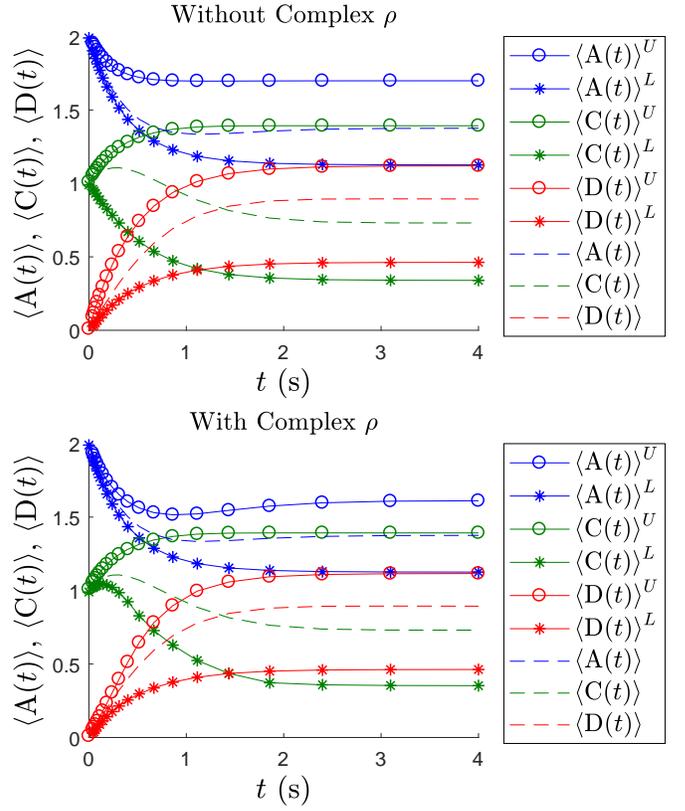}
    \caption{\label{Cyclic Bounds} Bounds on the mean molecular counts of species A, C, and D of Reaction System \eqref{cyclic system}.  Bounds on species B are omitted, because they are similar to those shown for species A.  The top panel shows bounds calculated without accounting for the imaginary components of the system\rq{}s eigenvalues, while the bottom panel shows the slight improvement that can be achieved by accounting for these imaginary components.  Both panels show the exact means calculated by directly solving the CME.}
\end{figure}

For example, consider the cyclic system 
\begin{equation} \label{cyclic system}
\begin{aligned}
\text{A} + \text{B} &\autorightarrow{$c_1$}{} \text{C} \\
\text{C} &\autorightarrow{$c_2$}{} \text{D} \\
\text{D} &\autorightarrow{$c_3$}{} \text{A} + \text{B}
\end{aligned}
\end{equation}
where the initial molecular counts are A $=2$, B $=1$, C $=1$, and D $=0$, and the rate constants are $c_1 = 1 \text{ s}^{-1}$, $c_2 = 1.1 \text{ s}^{-1}$, and $c_3 = 0.9 \text{ s}^{-1}$.
The smallest-magnitude eigenvalues of this system are $\lambda = 0, -2.1322 \pm 0.9741 i, -4.1637 \pm 1.5837i$.
If we follow the advice given in Section \ref{sec: choosing rho values}, and calculate bounds using $\mathcal R = \{0, -2.1322, -4.1637\}$, we obtain the bounds shown in the top panel of Figure \ref{Cyclic Bounds}.
However, by making use of the knowledge of the imaginary parts of the eigenvalues, we can produce the slightly improved bounds shown in the bottom panel.
The most notable improvements are for early times ($t < 0.5 \text{ s}$)

Given this potential to improve the bounds by using the imaginary parts of the low-magnitude eigenvalues, why has this paper been concerned almost solely with their real parts?
The fact is \lq\lq{}making use of the knowledge of the imaginary parts of the eigenvalues\rq\rq{} is not trivial.
One cannot simply use complex values of $\rho$ in SDPs \eqref{concrete optimization problem} and \eqref{concrete optimization problem, variance}.
The reason for this is that the argument for the derivation of LMIs \eqref{basic moment lmi, all t, integrated, integral inside} - \eqref{nonegative dependents,z} breaks down when $\rho$ is complex-valued.
It is possible to derive an analogous set of LMIs when $\rho$ is complex-valued, but this requires introducing entirely new classes of decision variables and constraints.
The resulting augmented versions of SDPs \eqref{concrete optimization problem} and \eqref{concrete optimization problem, variance} are considerably more complicated.
We felt that this extra complication would only distract from the main idea of this paper, and, as demonstrated by Figure \eqref{Cyclic Bounds}, it leads to only marginal improvement in the bounds.
Accordingly, we have deferred the discussion of how to account for complex eigenvalues to the supplementary material.

\section{Perfect Bounds in the Absence of the Closure Problem} \label{perfect bounds section} 
It is interesting to note that we can also apply our bounding method to stochastic chemical kinetic systems which do \emph{not} exhibit the closure problem, and that, doing so, it is possible to obtain perfect bounds.

\begin{figure}[H] 
    \centering
    \includegraphics[scale = 0.7]{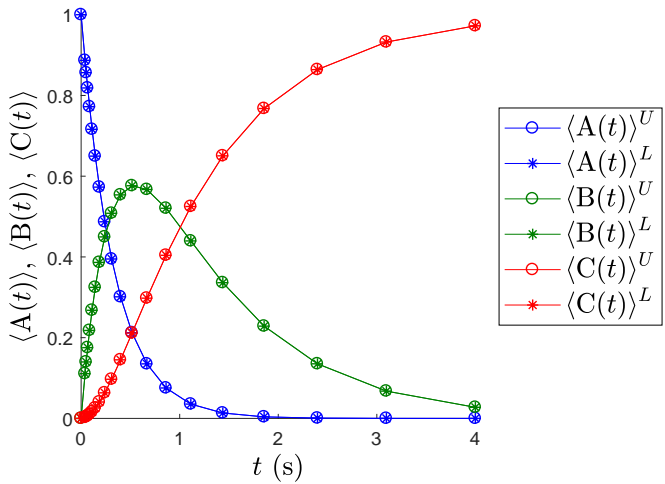}
    \caption{\label{collapsed bounds} Bounds on the mean molecular counts of species A, B, and C for Reaction System \eqref{unimolecular system}, which does not exhibit the closure problem.  For this example, the bounding method calculates perfect bounds, collapsing on the true mean trajectories.}
\end{figure}

For example, consider the reaction system
\begin{equation} \label{unimolecular system}
\text{A} \autorightarrow{$c_1$}{} \text{B} \autorightarrow{$c_2$}{} \text{C},
\end{equation}
where $c_1 = 1 \text{ s}^{-1}$, $c_2 = 3 \text{ s}^{-1}$, and there are initially $4$ molecules of A and $0$ molecules of each B and C.
Since every reaction in this system is unimolecular, it does not exhibit the closure problem.
The smallest-magnitude eigenvalues for this system are $\lambda = 0, -1, -3$.
Solving SDP \eqref{concrete optimization problem} and its minimization counterpart with $\mathcal R = \{ 0, -1, -3 \}$, we obtain the bounds shown in Figure \ref{collapsed bounds}.
The upper and lower bounding curves are indistinguishable from one another because there is essentially no gap between them; they have collapsed upon the true mean trajectories.

This is a rather nice feature of our bounding method, which, frankly, we did not expect.
We did not design the method with this collapsing behavior in mind.
However, as explained in the supplementary material, it naturally falls out of the math.
This example and others like it support the theoretical foundation of our bounding method -- in particular, the choice of exponential basis functions.

\section{Uncertainty in the Initial State} \label{uncertain initial state section}
In each of the foregoing examples, we have assumed that we knew the initial molecular count exactly.
This implies an initial probability distribution which is a Dirac distribution, where all the probability is concentrated on a single reachable state.
However, as suggested in Section \eqref{necessary steady-state conditions}, our method can also handle the more general situation where we don\rq{}t have exact knowledge of the initial molecular count, and the initial probability distribution (representing our knowledge of the system) is supported on several reachable states.
We demonstrate this capability with the following example.

Again, consider Reaction System \eqref{Ex6 system},
\begin{equation*} \label{Ex6 system, 2}
\text{A} + \text{B} \autorightarrow{$c_1$}{} \text{C} \autorightleftharpoons{$c_2$}{$c_3$} \text{D},
\end{equation*}
with the same rate constants given in Section \ref{mean bounding examples}.
In our prior analysis of this system, we 
assumed we knew the initial molecular counts A $=3$, B $=4$, C $=0$, and D $=0$.
This implies the set of reachable states $\mathcal X$ shown in Table \ref{reachable states table}.
\begin{table}
\centering
\begin{tabular}{c c}
\hline 
\hline
State 	& 	\ \ $\mbf x = (x_\text{A}, x_\text{B}, x_\text{C}, x_\text{D})$\ \\
\hline
1 		&	$(3, 4, 0, 0)$ 											\\
2 		&	$(2, 3, 1, 0)$ 											\\
3 		&	$(2, 3, 0, 1)$ 											\\
4		&	$(1, 2, 2, 0)$ 											\\
5		&	$(1, 2, 1, 1)$ 											\\
6		&	$(1, 2, 0, 2)$ 											\\
7		&	$(0, 1, 3, 0)$ 											\\
8		&	$(0, 1, 2, 1)$ 											\\
9		&	$(0, 1, 1, 2)$ 											\\
10		&	$(0, 1, 0, 3)$ 											\\
\hline
\hline
\end{tabular}
\caption{The set of reachable states $\mathcal X$ of the system described in Section \ref{mean bounding examples}.}
\label{reachable states table}
\end{table}
Furthermore, it implies an initial probability of zero for all states in Table \ref{reachable states table}, except State 1 which has an initial probability of one.



\begin{figure}[H] 
    \centering
    \includegraphics[scale = 0.7]{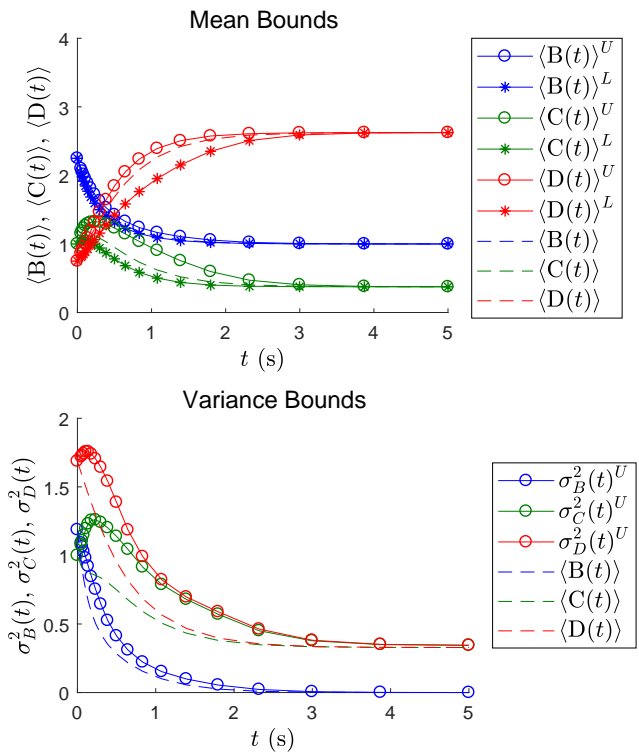}
    \caption{\label{distributed initial state} Bounds on the mean and variance for Reaction System \eqref{Ex6 system} with an uncertain initial state.}
\end{figure}

This time, we will assume uncertainty in the initial state, and we will express this uncertainty by assigning a nonzero initial probability to three distinct reachable states $\mbf x \in \mathcal X$.
In particular, we will assign initial probabilities of $\frac{1}{4}, \frac{1}{2}$, and $\frac{1}{4}$ to States 1, 4, and 10, respectively, with all other reachable states having an initial probability of zero.
Once we have decided on the set of species to be considered independent (e.g., species A and C), we can easily calculate the initial low-order moments $\bm \mu_L(0)$ corresponding to this initial distribution $P(\cdot, 0)$ using Equation \eqref{definition of moments}.
We can then apply SDPs \eqref{concrete optimization problem} and \eqref{concrete optimization problem, variance} to calculate bounds on the means and variances for this system over time.
For the sake of comparison to Figure \ref{fig: Ex6 Mean and Variance Bounds, more P}, we again use $m = 3$ and  $\mathcal R = \{ 0, -2, -2.4, -4.4 \}$.
The results are shown in Figure \ref{distributed initial state}.

The first thing to notice in comparing Figures \ref{fig: Ex6 Mean and Variance Bounds, more P} and \ref{distributed initial state} is that the starting point of each mean and variance trajectory is different between the two figures.
This is consistent with the fact that the initial distribution $P(\cdot, 0)$ and thus the initial moments $\mbf y_L(0)$ are different for the two figures.
The second thing to notice is that both plots approach the same steady-state at long times.
This is consistent with the fact that Reaction System \eqref{Ex6 system} has just one steady-state, in which species C and D are in equilibrium.
Finally, notice that the quality of the bounds is similar between the two plots. 
At least visually, the bounds in Figure \ref{distributed initial state} are just as tight as those in Figure \ref{fig: Ex6 Mean and Variance Bounds, more P}.
This may seem somewhat counter-intuitive given that Figure \ref{distributed initial state} was generated assuming uncertainty in the initial state.
However, recall that once this uncertainty is expressed in an initial probability distribution $P(\cdot, 0)$, the means and variances (which are \emph{expectation} values based on $P(\cdot,0)$) are precisely defined. 

\section{Conclusion}
This paper has described a method for calculating rigorous bounds on time-varying stochastic chemical kinetic systems.
In particular, we have formulated SDPs for calculating time-varying bounds on the mean molecular count of each species in the system and the variances in these counts.
This idea is an extension of the method described by several authors \cite{dowdy2017using, dowdy2018bounds, sakurai2017convex, kuntz2017rigorous, ghusinga2017exact} for calculating bounds on the steady-state (i.e., stationary) distribution of a stochastic chemical kinetic system.

As a proof of concept, we have demonstrated the bounding method for a simple stochastic chemical system for which analytical means and variances are available.
For this example, we have seen that our bounds are, in fact, valid.
Furthermore, we have seen that they can be very tight, given the appropriate choice of the parameter set $\mathcal R$.

We also applied the bounding method to a slightly more complicated reaction system, which demonstrates that method also applies to systems which reach a dynamic equilibrium at long times.
With this example, we saw that the bounds we obtain are not dramatically sensitive to the values of $\rho$ we select in our parameter set $\mathcal R$.

%

While the majority of the paper was written assuming that the parameter set $\mathcal R$ contained strictly real values $\rho$, in Section \ref{complex eigenvalue section} we saw that it is possible to obtain improved bounds by also using values of $\rho$ with nonzero imaginary parts -- though at the expense of solving a larger, more complicated SDP.

In Section \ref{perfect bounds section}, we saw an example which does not exhibit the closure problem, for which the bounds generated by our method collapse upon the true mean trajectories, supporting the theory underlying our approach.

Finally, in Section \ref{uncertain initial state section}, we demonstrated that our method can also handle the scenario when the initial state of the system is not known exactly and we instead have nonzero initial probabilities associated with several reachable states.

In theory, our bounding method could be applied to stochastic chemical kinetic systems of arbitrary size.
However, to do this, there are two practical issues that must be overcome:
first, we need to formalize a procedure for selecting the set $\mathcal R$;
second, we need to further explore options for mitigating the numerical issues mentioned in Section \ref{scaling section}.
Strategies for overcoming these issues will be the subject of a forthcoming publication.

Despite the method\rq{}s incompleteness, it is a theoretically novel, interesting approach to the closure problem in stochastic chemical kinetics.
We share it with the community in the hope that it might inspire further research in the area.

\section{Implementation Details}
All numerical examples in this paper were computed on a 64-bit Dell Precision T3610 workstation with a 3.70 GHz Intel Xeon CPU. 
In the example, CVX \cite{cvx} was used to model the SDP, using the default tolerance (i.e. “precision”) settings. 
SeDuMi \cite{S98guide} was used as the underlying solver.

\begin{acknowledgments}
Financial support from the Novartis-MIT Center for Continuous Manufacturing is gratefully acknowledged.
\end{acknowledgments}


\bibliography{myBib_dynamic}

\end{document}